\documentclass[a4paper, reqno,12pt]{amsart}
\usepackage{tikz, tikz-3dplot}
\usepackage{amssymb,latexsym, amsmath, amsthm, amscd, array, graphicx, color}
\usepackage{tikz}
\usepackage[english]{babel}
\usetikzlibrary{matrix,arrows,calc}
\usetikzlibrary{positioning}
\usepackage{parskip}
\usepackage[margin=1in]{geometry}
\usepackage{physics}
\usepackage[utf8]{inputenc}
\usepackage{amssymb}
\usepackage[export]{adjustbox}
\usepackage{setspace}
\usepackage{subcaption}
\usepackage{caption}
\usepackage{wrapfig}
\usepackage{bbm}
\usepackage{textcomp}
\usepackage{amsmath}
\usepackage{mathtools, xparse}
\usepackage{enumitem}
\usepackage{physics}
\usepackage{textcomp}
\usepackage{siunitx}
\usepackage{mathtools}
\usepackage{graphics}
\usepackage{mathrsfs}
\usepackage{url}
\theoremstyle{plain}
\newtheorem{theorem}{Theorem}
\newtheorem{lemma}[theorem]{Lemma}
\newtheorem{proposition}[theorem]{Proposition}

\theoremstyle{definition}
\newtheorem{remark}[theorem]{Remark}

\def\C{{\mathbb C}}
\def\Z{{\mathbb Z}}

\def\R{{\mathbb R}}
\def\N{{\mathbb N}}

\DeclarePairedDelimiterX{\inp}[2]{\langle}{\rangle}{#1, #2}
\newcommand{\unorm}[1]{{\left\vert\kern-0.25ex\left\vert\kern-0.25ex\left\vert #1 
    \right\vert\kern-0.25ex\right\vert\kern-0.25ex\right\vert}}

\title{Random Configurations and Fourier Decoupling}
\author{Siddharth Mulherkar}
	\address{Department of Mathematics, University of California, Los Angeles, USA}
	\email{sidmulherkar@g.ucla.edu}
\date{\today}
\begin{document}
\maketitle
\begin{abstract}
We use tools from high dimensional probability to analyze random configurations of balls in the context of Fourier decoupling and answer a question on the structure of superlevel sets in the setting of decoupling due to Guth. 
\end{abstract}

\maketitle
\section{Introduction}

In this article we will study random configurations in the context of Fourier decoupling theory and will answer a question of Guth \cite{GuthYouTubeSharpExamples2023}. Heuristics coming from probability theory often play an important guiding role in Fourier analysis. The notions of square root cancellation, pseudo-randomness and probabilistic models for various Fourier-analytic and number theoretic objects can often tell us how these objects \textit{should} behave. These heuristics can help us formulate conjectures and could also help us determine the types of inquiry that may yield progress on important problems. With this in mind, it can be fruitful to study random configurations in the context of Fourier analysis. For example, the work of Bourgain and Talagrand on $\Lambda(p)$ sets \cite{bourgain1989bounded} \cite{talagrand1994majorizing}, and the work of Mockenhaupt and Schlag \cite{mockenhauptSchlag2009majorantRandomSets} on the Hardy-Littlewood majorant problem have helped govern heuristics for the study of exponential sums and have also found applications in compressed sensing and signal processing (see \cite{IosevichMayeli2024LambdaSignalRecovery}, \cite{CandesPlan2010RISPless}, \cite{Tao2008UUP} for illustrative examples). We investigate ideas along these lines in the context of Fourier decoupling. A central question in this context is: what is the \textit{typical} behavior of natural Fourier analytic objects? In \cite{mulherkar2025random} the author investigated random configurations in the context of the Mizohata-Takeuchi conjecture, and showed that the Mizohata-Takeuchi conjecture is true \textit{most}\footnote{To be precise, it was shown that a generic random weight satisfies a sharp inequality of Mizohata-Takeuchi type. See \cite{mulherkar2025random}, \cite{cairo2025counterexample}, \cite{CairoZhang2025PowerLossMizohataTakeuchi} for context and discussion on this problem.} of the time. We use some of the tools developed in \cite{mulherkar2025random} and refine them to answer questions about the typical behavior of configurations that naturally arise in the context of Fourier decoupling.

\subsection{Setting of the Problem and Main Questions:}
We recall the setting of Fourier decoupling below. Let the set $\mathcal{N}_{1/N} (S^{1})$ be the ${1/N}$ neighborhood of the circle $S^{1} \subset \R^2$. Divide $\mathcal{N}_{\frac{1}{N}}(S^{1})$ into $N$ approximate rectangles of dimensions $1/N \times 1/N^2$ and call this collection $\{\theta\}$. Let $f \in \mathcal{S}(\mathbb{R}^2)$ (the Schwartz class on $\R^2$) be such that 
\begin{equation}\label{decouplingsetup}
    f = \sum_{\theta} f_\theta \quad \mathrm{with} \quad \mathrm{supp}(\hat{f_{\theta}}) \subseteq \theta,
\end{equation}
then a celebrated theorem of Bourgain and Demeter \cite{BourgainDemeter2015Decoupling} states that 
\begin{equation}
    \norm{f}_{L^6(\mathbb{R}^{2})} \lesssim_{\epsilon} N^{\epsilon} \left(\sum_{\theta} \norm{f_\theta}^{2}_{L^{6}(\mathbb{R}^2)}\right)^{1/2}. 
\end{equation}
Suppose in addition that for all $\theta$ we also have that
\begin{equation*}
     \norm{f_{\theta}}_{L^{\infty}(\R^2)} \leq 1.   
\end{equation*}
A consequence of more refined work by Guth, Maldague, and Wang on the decoupling problem yields a bound for the size of the super level sets $U_{\lambda}$. In this setting, define $U_{\lambda} :=\{x \in \R^2: \abs{f(x)}> \lambda\}$. Then if $B_{N^2} \coloneqq \{x: \abs{x} \leq N^2\}$ (that is the ball of radius $N^2$ in $\R^2$), then in \cite{GuthMaldagueWang2020ImprovedParabola} it was shown that
\begin{equation}\label{superlevel}
    \abs{U_{\lambda} \cap B_{N^2}} \lesssim_{\epsilon} N^{\epsilon} \frac{N^7}{\lambda^{6}}.
\end{equation}

In light of this bound we are led to the following natural questions:

\textbf{Questions:}
\begin{enumerate}
    \item For different values of $\lambda$, what kinds of functions $f$ saturate the bound (\ref{superlevel})? 
    \item What can we say about the structure of $U_{\lambda}$ if the bound (\ref{superlevel}) is saturated?
\end{enumerate}
 One can think of these questions as asking for an `inverse theorem' in the context of Fourier decoupling. Inverse theorems have played a significant role in recent developments in additive combinatorics and number theory (see this excellent talk of Gowers \cite{GowersYouTubeInverseTheoremsAdditiveCombinatorics} for more background), and a satisfactory answer to these questions could lead to important insights about the behavior of exponential sums and related questions (see this survey \cite{Guth2022DecouplingSurvey} of Guth for more discussion). 

When $\lambda = N^{1/2}$ it is easy to find many examples where the bound (\ref{superlevel}) is saturated. On the other hand if $N^{1/2} < \lambda \leq N$, there are very few examples of functions that saturate the bounds in (\ref{superlevel}), and it is not clear whether this bound is sharp for all values of $\lambda$. We remark that the question of the sharpness of (\ref{superlevel}) was investigated in \cite{FuGuthMaldague2023}, but in the case where $S^{1}$ was replaced by the parabola $\mathbb{P}^{1}$, where the authors were able to leverage techniques coming from the study of some well known exponential sums. These techniques do not carry over to the case of the circle $S^{1}$. In \cite{GuthYouTubeSharpExamples2023} Guth asked whether the bound (\ref{superlevel}) can be saturated by a random set. In particular he asked the following:
\begin{enumerate}
    \item If $V$ is a union of $N^{7}\lambda^{-6}$ random unit balls, and $\lambda > N^{1/2}$, is it typically possible (in the setting above) that $\norm{f}_{L^{1}(V)} \sim \lambda \abs{V}$?
\end{enumerate}
If the answer to this question was yes, then it would be that in the case when (\ref{superlevel}) was saturated, the set $U_{\lambda} \cap B_{N^2}$ could look like a random set. However, we answer this question in the negative. In fact we show that on the random set $V$,  $f$ will be strictly smaller than what the superlevel set bound (\ref{superlevel}) guarantees. To set up our main theorem we first define the class of functions $\mathcal{F}$ as follows:
\begin{equation}
    \mathcal{F} \coloneqq \left\{f \in \mathcal{S}(\R^2): f= \sum_{\theta} f_{\theta}, \quad \mathrm{supp} (\hat{f_{\theta}}) \subseteq \theta, \quad \norm{f_{\theta}}_{L^{\infty}(\R^2)} \leq 1 \right\},
\end{equation}
where the collection $\{\theta\}$ is defined as in (\ref{decouplingsetup}). Our main theorem gives a precise asymptotic for $\sup_{f \in \mathcal{F}}\int_{V}\abs{f(x)}dx$. 
\begin{theorem}[Corollary of Main Theorem, Informal Statement]
    In the setting above, let $V_\lambda$ be a random collection of $N^{7}\lambda^{-6}$ unit balls contained in $B_{N^2}$, then with high probability 
    \begin{equation}\label{maintheoremcases}
        \sup_{f \in \mathcal{F}}\norm{f}_{L^1(V_{\lambda})} \sim_{\epsilon} 
        \begin{cases}
            N^{1/2}\abs{V_{\lambda}} \quad \text{when } N^{1/2} < \lambda \leq N^{5/6} \\
            \frac{\lambda^3}{N^2} \abs{V_{\lambda}} \quad \text{when } N^{5/6} \leq \lambda \leq N
        \end{cases}
    \end{equation}
\end{theorem}
 In any case one can deduce that whenever $N^{1/2} \leq \lambda < N$, that
\begin{equation}
    \sup_{f \in \mathcal{F}}\norm{f}_{L^1(V_\lambda)} \sim o(\lambda) \abs{V_\lambda}.
\end{equation}
What this shows is that a set $U_{\lambda}$ that saturates (\ref{superlevel}) cannot look random. From here, one can deduce that examples of sets $U_{\lambda} \cap B_{N^2}$ that saturate (\ref{superlevel}) must possess some kind of special $\textit{structure}$. Finding examples of such sets and determining what kind of structure $U_{\lambda} \cap B_{N^2}$ must possess remains an important challenge. See \cite{Guth2022DecouplingSurvey} for related discussion.

The two different regimes for the estimate (\ref{maintheoremcases}) may appear mysterious to the reader. We clarify this by proving a more general estimate for collections of random sets of arbitrary cardinality. 

\begin{theorem}[Main Theorem, Informal Statement]\label{MAINtheoreminformal}
    For each $N \leq q \leq CN^4$, let $V$ be a random collection of $q$ unit balls contained inside $B_{N^2}$, then with high probability
    \begin{equation}\label{maintheoremeq}
        \sup_{f \in \mathcal{F}}\norm{f}_{L^1(V)} \sim_{\epsilon} N^{\epsilon}(N^{1/2}q + N^{3/2}\sqrt{q}). 
    \end{equation}
\end{theorem}
    when $q = N^{7}\lambda^{-6}$, one can check that the first term in (\ref{maintheoremeq}) dominates whenever $N^{1/2} \leq \lambda \leq N^{5/6}$ and then second term in (\ref{maintheoremeq}) dominates when $N^{5/6} < \lambda \leq N$. These observations yield the estimate (\ref{maintheoremcases}).

\subsection{Setup for Theorem \ref{MAINtheoreminformal}}\label{setup}

In this section we will convert the informal statement of Theorem \ref{MAINtheoreminformal} into a formal statement. Decompose $B_{N^2}$ into a finitely overlapping union of $CN^4$ unit balls $\{V_k\}_{k=1}^{CN^{4}}$ for some absolute constant $C$. For each $1 \leq k \leq CN^{4}$, define the Bernoulli random variable $\delta_k$ as follows:
\begin{equation}\label{selector}
    \mathbb{P}(\delta_k=1) = O\left(\frac{q}{N^4}\right) \coloneqq \delta, \quad \mathbb{P}(\delta_{k} =0) = 1- \delta.
\end{equation}
Then with high probability, the random set $V$ defined as
\begin{equation}\label{randomsetballs}
    V \coloneqq \left\{ \bigcup_{k: \delta_k =1} V_{k} \right\} \subseteq B_{N^2},
\end{equation}
is such that $\abs{V} \sim q$. 
\begin{theorem}\label{MAINtheorem}
    In the setting coming from the discussion above,
    \begin{equation}\label{mainproptechnical}
        \mathbb{E}\sup_{f \in \mathcal{F}} \norm{f}_{L^1(V)}= \mathbb{E}\sup_{f \in \mathcal{F}}\sum_{k=1}^{CN^4} \delta_k \int_{V_k} \abs{f(x)}dx \sim_{\epsilon} N^{\epsilon}(N^{1/2}q + N^{3/2} \sqrt{q}).
    \end{equation}
\end{theorem}
The first equality can easily be seen to be true by chasing definitions, the second inequality is the main technical result of this article.

The lower bound in (\ref{mainproptechnical}) comes from some probabilistic constructions involving wave packets, which we show in Section \ref{lowerboundsection}. The technical heart of this paper, however, is in proving the upper bound which we restate here.
\begin{proposition}\label{mainprop}
We have the following inequality:
    \begin{equation}\label{mainpropeq}
        \mathbb{E}\sup_{f \in \mathcal{F}}\sum_{k=1}^{CN^4} \delta_k \int_{V_k} \abs{f(x)}dx \lesssim_{\epsilon} N^{\epsilon}(N^{1/2}q + N^{3/2}\sqrt{q})
    \end{equation}
\end{proposition}
\begin{remark}
    One could ask about the case when $1 \leq q \leq N$, however, this regime is not so interesting since in this case the random set $V$ is quite sparse and there exist $f \in \mathcal{F}$ for which $f$ is large on all of $V$. In fact we have that when $1 \leq q \leq N$
    \begin{equation*}
        \mathbb{E}\sup_{f \in \mathcal{F}}\norm{f}_{L^1(V)} \sim q N.
    \end{equation*}
    The upper bound is trivial since $\norm{f}_{L^{\infty}} \leq N$ for all $f \in \mathcal{F}$. Moreover a `wave packet tuning' argument of Guth \cite{GuthYouTubeSharpExamples2023}, produces a corresponding lower bound.
\end{remark}

\subsection{Suprema of Stochastic Processes:}

One of the key inspirations behind the proof of Proposition \ref{mainprop} is to use ideas from the theory of bounding the expected supremum of stochastic processes. We can reframe (\ref{mainpropeq}) more abstractly as a question that asks for tight bounds for the quantity
\begin{equation}\label{expectedsup}
    \mathbb{E}\sup_{t \in T} X_t,
\end{equation}
for some stochastic process $(X_t)_{t \in T}$. There is an incredibly rich and interesting theory surrounding this family of problems and we refer the reader to Talagrand's book \cite{talagrand2014upper} and the references therein for a comprehensive introduction. The study of the expected suprema of Gaussian processes is especially comprehensive; when $(X_t)_{t \in T}$ is a Gaussian process it is known that the quantity (\ref{expectedsup}) is completely determined by the metric induced by the $L^2$ distance between random variables on the index set $T$, given by $d(s,t) = (\mathbb{E}(X_t-X_s)^2)^{1/2}$ for $s,t \in T$. A powerful strategy to used bound the expected supremum of a Gaussian process is called the \textit{chaining} method. An informal description of chaining for Gaussian processes is as follows: assuming the process is centered\footnote{That is, $\mathbb{E}X_{t} =0$ for all $t \in T$.} we can write
\begin{equation*}
    \mathbb{E}\sup_{t \in T}X_t = \mathbb{E}\sup_{t \in T} (X_t-X_{t_0}),
\end{equation*}
and then write $(X_t-X_{t_0}) = \sum_{k} (X_{\pi_{k}(t)}-X_{\pi_{k-1}(t)}) $ as some telescoping sum, where the $\pi_{k}(t)$'s can be thought of as successive approximations to $t$ in the metric space $(T,d)$. One can then use measure concentration bounds\footnote{In the case of a Gaussian process, one usually uses the Gaussian tail bound: $\mathbb{P}(\abs{X_s-X_t}> \lambda) \leq 2 \exp\left(\frac{-\lambda^2}{2d(s,t)^2}\right)$.} to control the size of the increments $X_{\pi_{k}(t)}- X_{\pi_{k-1}(t)}.$  A celebrated application of the chaining method is \textit{Dudley's inequality} which states in this setting that
\begin{equation*}
    \mathbb{E} \sup_{t \in T}X_t \lesssim \int_{0}^{\infty} \sqrt{\log N(T,d, \epsilon)} d \epsilon,
\end{equation*}
where $N(T,d, \epsilon)$ is the $\epsilon$ - \textit{covering number} of $T$ (sometimes referred to as the metric entropy of $T$), that is, the minimum number of balls of radius $\epsilon$ in the metric space $(T,d)$ required to cover $T$. We refer the reader to the references \cite{talagrand2014upper} and \cite{vershynin2018high} for excellent introductions to the chaining method and applications. Our method crucially relies on Dudley's inequality, but it cannot be used directly. We need to manipulate the stochastic process we are interested in, and then compute entropy estimates in an appropriate way. We refer the reader to Appendix B for some more background on the precise statements from probability theory that we use in this paper.

\begin{remark}
    The astute reader may have realized that the quantity
    \begin{equation*}
        \mathbb{E} \sup_{t \in T}X_t
    \end{equation*}
    may not be well defined if $T$ is uncountable (since the random variable $\sup_{t \in T}X_t$ may not even be measurable). The fix for this is to alter the definition of the expected supremum of a stochastic process so that all the previous arguments make sense by defining
    \begin{equation*}
        \mathbb{E}\sup_{t \in T}X_t \coloneqq \sup_{F \subseteq T, \text{ finite}} \mathbb{E}\sup_{t \in F}X_t.
    \end{equation*}
    We will proceed with this convention for this paper as well.
\end{remark}

\subsection{Structure and Organization of the Paper:} 
\begin{itemize}
    \item In Section \ref{reductionsection}, we reduce the proof of Proposition \ref{mainprop} to another proposition, Proposition \ref{mainmainpropmodified}.
    \item  In Section \ref{upperboundsection}, we prove Proposition \ref{mainmainpropmodified}, which completes the proof of the upper bound for Theorem \ref{MAINtheorem}.
    \item In Section \ref{lowerboundsection}, we supply the proof of the lower bound for Theorem \ref{MAINtheorem}, using wave packets.
    \item In Section \ref{furtherremarkssection}, we give some further remarks on the paper and the proof strategy. This includes some suggestions for reading this paper. Once familiar with the basic setup of this problem, the reader is recommended to look at Section \ref{furtherremarkssection} for suggestions on how to read this paper.
    \item The appendix contains two parts. In Appendix A, we prove some technical propositions in Fourier analysis. In Appendix B, we give some background on the machinery coming from tools in high dimensional probability we use in the proof.
\end{itemize}

\subsection{Notation} We will use the following notation in this paper:
\begin{itemize}
    \item $A \lesssim B$ if there exists an absolute constant $C>0$ such that $A \leq C B$.
    \item $A \sim B$ if there exists an absolute constant such that $\frac{1}{C}B \leq A \leq CB$.
    \item $A \lesssim_{\epsilon} B$ if for every $\epsilon >0 $ there exists a constant $C_{\epsilon}>0$ (possibly depending on $\epsilon$) such that $A \leq C_{\epsilon} B$ for all $\epsilon >0$. For example, $\log R \lesssim_{\epsilon} R^{\epsilon}$ for all $R>0$.
    \item Given any $A, B \subseteq \mathbb{R}^n$, define $\mathcal{N}(A, B)$ to be the minimum number of translates of $B$ needed to cover $A$. If $\norm{\cdot}$ is a norm on $\mathbb{R}^n$, we write $\mathcal{N}(A, \norm{\cdot}, \epsilon)$ to be the minimum number of ball of radius $\epsilon$ in $\norm{\cdot}$ needed to cover $A$.  Covering numbers of the form $\mathcal{N}(A,B)$ are a basic tool in the study of suprema of stochastic processes and will be used frequently in this work.
\end{itemize}
\subsection{Tool and Computational Resource Disclosure} All ideas, theorems and technical statements in this paper are due to author, and were written by the author. Once a first draft of the paper was written, we used a large language model to review the paper. The author incorporated several suggestions from the language model, which included corrections of expository and mathematical errors, as well as minor simplifications to some arguments. 
\subsection{Acknowledgments:} The author would like to thank Terence Tao and Hong Wang for feedback and support. Part of this work was completed when the author was visiting the Institut Hautes \'{E}tudes des Scientifiques, and we thank the IHES for excellent working conditions. 
\section{Reduction to Proposition \ref{mainmainpropmodified}}\label{reductionsection}
The goal of this section is to reduce the statement of Proposition \ref{mainprop} to a similar statement, Proposition \ref{mainmainpropmodified} which we state below. 

For $j \in \N$ and each multi-index $\alpha = (i_1,i_2,...,i_j)$ where each $i_m \in \{1,2\}$, define $\abs{\alpha} = j$. For any $f \in \mathcal{S}(\R^2)$, define $D^{(\alpha)}f$ as follows:
\begin{equation*}
    D^{(\alpha)}f \coloneqq \frac{\partial f}{\partial x_{i_1} \partial x_{i_2} ... \partial x_{i_j}}.
\end{equation*}
  Recall from that $\{V_{k}\}$ is a collection of $CN^{4}$ unit balls. For each $k$, decompose $V_k$ into $\sim N^{\epsilon}$ many finitely overlapping unit balls ${V_{k,l}}$ of radius $N^{-\epsilon}$. Write $V_{k,l} =B(y_{kl}, N^{-\epsilon})$ as the ball of radius $N^{-\epsilon}$ around $y_{kl}$. For each $k$, we can label the $y_{kl}$ such that $1 \leq l \leq c_1N^{\epsilon}$, for some absolute constant $c_1$.
\begin{proposition}\label{mainmainpropmodified}
    Fix $\epsilon >0$ and let $(\delta_{k})$ be the collection of selector random variables defined in (\ref{selector}). Recall that
    \begin{equation}
        \mathbb{P}(\delta_k = 1) \sim \frac{q}{N^4}.
    \end{equation}
    Then for each multi-index $\alpha$, with $\abs{\alpha} \lesssim \frac{1}{\epsilon}$, and each $l$, we have that 
    \begin{equation}\label{newprop}
        \mathbb{E} \sup_{f \in \mathcal{F}} \sum_{k=1}^{CN^{4}} \delta_k \abs{D^{(\alpha)}f(y_{kl})} \lesssim_{\epsilon} N^{\epsilon} (N^{1/2}q + N^{3/2} \sqrt{q}).  
    \end{equation}
\end{proposition}
\begin{proof}[Proof that Proposition \ref{mainmainpropmodified} implies Proposition \ref{mainprop}]
    Assuming (\ref{newprop}) is true, we need to show that 
    \begin{equation}\label{reducedmainprop}
        \mathbb{E}\sup_{f \in \mathcal{F}}\sum_{k=1}^{CN^4} \delta_k \int_{V_k} \abs{f(x)}dx \lesssim_{\epsilon} N^{\epsilon} (N^{1/2}q + N^{3/2} \sqrt{q}).
    \end{equation}
    For each $k$, decompose $V_{k}$ into the finite overlapping collection of balls $\{V_{k,l}\}$, so that 
    \begin{equation}\label{subadditive}
        \int_{V_{k}} \abs{f(x)}dx \lesssim \sum_{l} \int_{V_{k,l}} \abs{f(x)}dx.
    \end{equation}
    For each $x \in V_{k,l}$, note that $x \in B(y_{kl}, N^{-\epsilon})$, for each $f \in \mathcal{F},$ we have by Taylor's theorem that
    \begin{equation}\label{firsttaylor}
        \abs{f(x)} \lesssim \left(\sum_{\abs{\alpha} \leq M} \abs{D^{(\alpha)}f(y_{kl})} \frac{N^{-\abs{\alpha} \epsilon}}{\alpha !} \right) +  \sup_{\abs{\alpha} = M+1} \abs{D^{(\alpha)}f(y_{kl})} \cdot O(N^{-M \epsilon}),
    \end{equation}
    for some $M \in \mathbb{N}$ (to be chosen later)
    The next lemma shows that we can bound the remainder uniformly for all $f \in \mathcal{F}$.
    \begin{lemma}
        For each $f \in \mathcal{F}$, we have
        \begin{equation}
            \sup_{f \in \mathcal{F}} \abs{D^{(\alpha)}f(y_{kl})} \lesssim_{\abs{\alpha}} N.
        \end{equation}
    \end{lemma}
    \begin{proof}
        Indeed, by Bernstein's inequality \cite[Chapter 5]{Wolff2003}, since each $f \in \mathcal{F}$ is such that $\hat{f}$ has compact support, we have
        \begin{equation*}
            \norm{D^{(\alpha)}f}_{\infty} \lesssim_{\abs{\alpha}} \norm{f}_{\infty} 
        \end{equation*}
    Moroever for each $f \in \mathcal{F}$,
    \begin{equation*}
        \norm{f}_{\infty} = \norm{\sum_{\theta}f_{\theta}}_{\infty} \leq \sum_{\theta}\norm{f_{\theta}}_{\infty} \leq N,
    \end{equation*}
    which proves the claim.
    \end{proof}
    In light of the previous lemma, for each $f \in \mathcal{F}$,  (\ref{firsttaylor}) can be rewritten as
    \begin{align}
        \abs{f(x)} \lesssim \sum_{\abs{\alpha} \leq M} \abs{D^{(\alpha)}f(y_{kl})} \frac{N^{-\abs{\alpha} \epsilon}}{\alpha !} + O_{M}(N^{-M \epsilon +1})
    \end{align}
    Set $M = \lceil \frac{10001}{\epsilon}\rceil$, so that the above expression can be rewritten as
    \begin{equation}\label{finaltaylor}
        \abs{f(x)} \lesssim \sum_{\abs{\alpha} \lesssim \frac{1}{\epsilon} } \abs{D^{(\alpha)}f(y_{kl})} \frac{N^{-\abs{\alpha} \epsilon}}{\alpha !} + O_{\epsilon}(N^{-10000}).
    \end{equation}
    We now prove (\ref{reducedmainprop}) assuming (\ref{newprop}). First, by (\ref{subadditive}) we have that
    \begin{equation}
        \mathbb{E}\sup_{f \in \mathcal{F}}\sum_{k=1}^{CN^4} \delta_k \int_{V_k} \abs{f(x)}dx \lesssim \mathbb{E}\sup_{f \in \mathcal{F}}\sum_{k=1}^{CN^4}  \delta_k \sum_{l} \int_{V_{k,l}} \abs{f(x)}dx. 
    \end{equation}
    Now, applying the estimate (\ref{finaltaylor}) yields 
    \begin{align*}
        &\mathbb{E}\sup_{f \in \mathcal{F}}\sum_{k=1}^{CN^4} \delta_k \sum_{l} \int_{V_{k,l}} \abs{f(x)}dx \\ 
        &\lesssim \mathbb{E}\sup_{f \in \mathcal{F}}\sum_{k=1}^{CN^4} \delta_k \sum_{l} \int_{V_{k,l}} \left(\sum_{\abs{\alpha} \lesssim \frac{1}{\epsilon} } \abs{D^{(\alpha)}f(y_{kl})} \frac{N^{-\abs{\alpha} \epsilon}}{\alpha !}  + O_{\epsilon}(N^{-10000}) \right)dx . \\
        &\lesssim  \sum_{l} \sum_{ \abs{\alpha} \lesssim \frac{1}{\epsilon}}  \mathbb{E} \sup_{f \in \mathcal{F}} \left( \sum_{k=1}^{CN^4}   \delta_k \frac{N^{-\abs{\alpha}\epsilon}}{\alpha !}\int_{V_{k,l}}\abs{D^{(\alpha)} f(y_{kl})} dx \right) + O_{\epsilon} (N^{-1000+\epsilon})\\
    \end{align*}
    where in the second line, we have used the estimate (\ref{finaltaylor}), and in the third line, we have interchanged sums, and pulled the sums in $l$ and $\alpha$ outside the supremum over $f \in \mathcal{F}$. Noting that for each $k$ and $l$, the area of $V_{k,l}$ is $\sim N^{-\epsilon}$, we get that the above expression is bounded by
    \begin{equation*}
        \lesssim  \sum_{l} \sum_{ \abs{\alpha} \lesssim \frac{1}{\epsilon}} \mathbb{E} \sup_{f \in \mathcal{F}} \left( \sum_{k=1}^{CN^4}   \delta_k \frac{N^{-(\abs{\alpha}+1)\epsilon}}{\alpha !}\abs{D^{(\alpha)} f(y_{kl})} \right) + O_{\epsilon}(N^{-1000+\epsilon}).
    \end{equation*}
    Now noting that the term $\frac{N^{-(\abs{\alpha}+1)\epsilon}}{\alpha !}$ for each $\alpha$, we see that the above expression is bounded by 
    \begin{equation}\label{almosttherelemma}
        \leq  \sum_{l=1}^{c_1N^{\epsilon}} \sum_{ \abs{\alpha} \lesssim \frac{1}{\epsilon}} \mathbb{E} \sup_{f \in \mathcal{F}} \left( \sum_{k=1}^{CN^4}   \delta_k\abs{D^{(\alpha)} f(y_{kl})} \right) + O_{\epsilon} (N^{-1000+\epsilon}).
    \end{equation}
    Now we apply (\ref{newprop}), which gives us an estimate for the inner term.
    \begin{equation*}
        \mathbb{E} \sup_{f \in \mathcal{F}} \left( \sum_{k=1}^{CN^4}   \delta_k\abs{D^{(\alpha)} f(y_{kl})} \right) \lesssim_{\epsilon} N^{\epsilon}(N^{1/2}q + N^{3/2}\sqrt{q})
    \end{equation*}
    Finally, noting that the number summations in (\ref{almosttherelemma}) in $l$ and $\abs{\alpha}$ combined is $\lesssim_{\epsilon} N^{\epsilon}$, this yields that the expression (\ref{almosttherelemma}) is $\lesssim_{\epsilon} N^{\epsilon}(N^{1/2}q + N^{3/2}\sqrt{q})$, so that 
    \begin{equation*}
        \mathbb{E}\sup_{f \in \mathcal{F}}\sum_{k=1}^{CR^d} \delta_k \int_{V_k} \abs{f(x)}dx \lesssim_{\epsilon} N^{\epsilon}(N^{1/2}q + N^{3/2}\sqrt{q}).
    \end{equation*}
    which proves the claim.
\end{proof}
\section{Proof of the Upper Bound}\label{upperboundsection}
\subsection{Symmetrization and Further Reductions}
        The previous section shows that we need to establish the desired bound
        \begin{equation}\label{nonsymmetrized}
            \mathbb{E}\sup_{f \in \mathcal{F}} \sum_{k=1}^{CN^4} \delta_k \abs{D^{(\alpha)}f(y_{kl})}  \lesssim_{\epsilon} N^{\epsilon}(N^{1/2}q + N^{3/2}\sqrt{q})
        \end{equation}
        where 
        \begin{equation}
            q \in \N, \quad \text{and} \quad N \leq q \leq CN^4, 
        \end{equation}
        where $q$ is the number of random points chosen from $\{y_{kl}\}_{k=1}^{CN^4}$, in expectation. The next step is a step known as `symmetrization' in the literature (see \cite[Lemma 6.3.2]{vershynin2018high}). 
        First we write
        \begin{equation*}
            \delta_k = (\delta_k - \delta) + \delta
        \end{equation*}
        Inserting this into (\ref{nonsymmetrized}) we have 
        \begin{equation*}
            \mathbb{E}\sup_{f \in \mathcal{F}} \sum_{k=1}^{CN^4} \delta_k\abs{D^{(\alpha)}f(y_{kl})} \leq \mathbb{E}\sup_{f \in \mathcal{F}} \abs{ \sum_{k=1}^{CN^4}(\delta_k- \delta)\abs{D^{(\alpha)}f(y_{kl})}} + \delta \sup_{f \in \mathcal{F}} \sum_{k=1}^{CN^4} \abs{D^{(\alpha)}f(y_{kl})}. 
        \end{equation*}
        The second term is deterministic and easy to estimate. First note that 
        \begin{equation}\label{deltabasicbound}
            \delta \sim \frac{q}{N^4},
        \end{equation}
        and some basic Fourier analytic arguments give that
        \begin{equation}\label{localconstl1}
            \sup_{f \in \mathcal{F}} \sum_{k=1}^{CN^4} \abs{D^{(\alpha)}f(y_{kl})} \lesssim \norm{D^{(\alpha)}f}_{L^1\left(B_{N^2}\right)} \lesssim \norm{f}_{L^1\left(B_{N^2}\right)} \lesssim_{\epsilon} N^{9/2 + \epsilon}. 
        \end{equation}
        We prove these estimates rigourously in the appendix, but we can explain the basic heuristic here. Since every $f \in \mathcal{F}$ has Fourier transform supported in $[-2,2]^2$, every $f \in \mathcal{F}$ is `locally constant'. Thus one can make the approximation $\sum_{k}\abs{D^{(\alpha)}f(y_{kl})} \lesssim \norm{D^{(\alpha)}f}_{L^1\left(B_{N^2}\right)}$. The third inequality follows from Bernstein's inequality and the final inequality comes from integrating in the superlevel set bound (\ref{superlevel}). Combining the estimates (\ref{deltabasicbound}) and (\ref{localconstl1}) we obtain
        \begin{equation}
            \delta \sup_{f \in \mathcal{F}} \sum_{k=1}^{CN^4} \abs{D^{(\alpha)}f(y_{kl})} \lesssim_{\epsilon} N^{1/2+ \epsilon} q,
        \end{equation}
        which corresponds to the first term on the right hand side of (\ref{nonsymmetrized}). It remains to establish
        \begin{equation}
            \mathbb{E}\sup_{f \in \mathcal{F}} \abs{\sum_{k=1}^{CN^4} (\delta_k - \delta) \abs{D^{(\alpha)}f(y_{kl})}} \lesssim N^{3/2 + \epsilon} \sqrt{q}
        \end{equation}
        Let $(\delta_k')$ be collection of selector random variables with mean $\delta$ which is independent from the collection $(\delta_k)$. We can write 
        \begin{align*}
            \mathbb{E}\sup_{f \in \mathcal{F}} \abs{\sum_{k=1}^{CN^4} (\delta_k - \delta) \abs{D^{(\alpha)}f(y_{kl})}} &= \mathbb{E}_{\delta}\sup_{f \in \mathcal{F}} \abs{\mathbb{E}_{\delta'} \sum_{k=1}^{CN^4} (\delta_k- \delta_k')\abs{D^{(\alpha)}f(y_{kl})}}\\
            &\leq \mathbb{E}_{\delta, \delta'} \sup_{f \in \mathcal{F}} \abs{\sum_{k=1}^{CN^4}(\delta_k - \delta_k') \abs{D^{(\alpha)}f(y_{kl})}} 
        \end{align*}
        where the last inequality follows  from the fact that $\sup_{t\in T}\mathbb{E}X_t \leq \mathbb{E}\sup_{t \in T}X_{t}$, and the independence of the random variables $(\delta_k)_{k=1}^{n}$ and $(\delta_k')_{k=1}^{n}$.
        Let $(\epsilon_{k})$ be an independent and identically distributed collection of Rademacher random variables, that is $\mathbb{P}(\epsilon_k = \pm 1) = \frac{1}{2}$, and note that the random variables 
        \begin{equation*}
            \delta_k - \delta_k' \quad \text{and} \quad \epsilon_k(\delta_k- \delta_k')
        \end{equation*}
        have the same law. Therefore,
        \begin{align*}
            &\mathbb{E}_{\delta, \delta'} \sup_{f \in \mathcal{F}} \abs{\sum_{k=1}^{CN^4}(\delta_k - \delta_k') \abs{D^{(\alpha)}f(y_{kl})}} = \mathbb{E}_{\delta, \delta', \epsilon} \sup_{f \in \mathcal{F}} \abs{\sum_{k=1}^{CN^4} \epsilon_k(\delta_k - \delta_k') \abs{D^{(\alpha)}f(y_{kl})}}\\
        \end{align*}
        and the previous term can be bounded as follows:
        \begin{align}\label{symmetrized}
            \mathbb{E}_{\delta, \delta', \epsilon} \sup_{f \in \mathcal{F}} \abs{\sum_{k=1}^{CN^4} \epsilon_k(\delta_k - \delta_k') \abs{D^{(\alpha)}f(y_{kl})}} &\leq \mathbb{E} \sup_{f \in \mathcal{F}} \abs{\sum_{k=1}^{CN^4} \epsilon_k\delta_k \abs{D^{(\alpha)}f(y_{kl})}} + \mathbb{E} \sup_{f \in \mathcal{F}} \abs{\sum_{k=1}^{CN^4} \epsilon_k\delta_k \abs{D^{(\alpha)}f(y_{kl})}} \notag\\
            &\leq 2 \mathbb{E} \sup_{f \in \mathcal{F}}\abs{\sum_{k=1}^{CN^4} \epsilon_k\delta_k \abs{D^{(\alpha)}f(y_{kl})}}. 
        \end{align}
        Let $S$ be the (random) set defined as
        \begin{equation}
            S \coloneqq \{k: \delta_k = 1 \}.
        \end{equation}
        The next step (which also appears several times in the literature, see for example, \cite[Theorem 8.3.15]{vershynin2018high}) is to condition on the random set $S$. By Fubini's theorem, the expression (\ref{symmetrized}) can be written as
        \begin{equation}\label{conditional}
            \mathbb{E}_{\delta} \mathbb{E}_{\epsilon} \sup_{f \in \mathcal{F}} \abs{\sum_{k \in S} \epsilon_k\abs{D^{(\alpha)}f(y_{kl})}}
        \end{equation}
        The main result of the next section is the following estimate:
        \begin{equation}\label{symmetrizedprocessmain}
            \mathbb{E}_{\epsilon} \sup_{f \in \mathcal{F}}\abs{\sum_{k \in S} \epsilon_k  \abs{D^{(\alpha)}f(y_{kl})}} \lesssim_{\epsilon} N^{3/2+ \epsilon} \sqrt{\abs{S}},
        \end{equation}
            where this estimate holds for \textit{any} subset $S \subseteq [CN^4]$. Plugging this estimate into (\ref{conditional}), we have
        \begin{equation*}
            \mathbb{E}_{\delta} \mathbb{E}_{\epsilon} \sup_{f \in \mathcal{F}} \abs{\sum_{k \in S} \epsilon_k \abs{D^{(\alpha)}f(y_{kl})}} \lesssim_{\epsilon} \mathbb{E}_{\delta}  (N^{3/2+ \epsilon} \sqrt{\abs{S}}) \leq N^{3/2+ \epsilon} \sqrt{q}
        \end{equation*}
        where the final inequality follows from Jensen's inequality. This establishes the required upper bound, conditional on proving the estimate (\ref{symmetrizedprocessmain}), which we will prove in the next section.
        \subsection{Bounding the Symmetrized Process}
        In the next two sections we establish the estimate (\ref{symmetrizedprocessmain}), which we recall here.
        \begin{proposition}
            Fix any subset $S \subseteq [CN^4]$. We have the following estimate:
            \begin{equation*}
                \mathbb{E}_{\epsilon} \sup_{f \in \mathcal{F}}\abs{\sum_{k \in S} \epsilon_k  \abs{D^{(\alpha)}f(y_{kl})}} \lesssim_{\epsilon} N^{3/2+ \epsilon} \sqrt{\abs{S}}.
            \end{equation*}
            \end{proposition}
            We begin by making a reduction, where we reduce to studying a similar problem where the $\sup$ is taken over a simpler class $\mathcal{F_{\theta}}$, defined as
            \begin{equation*}
                \mathcal{F}_{\theta} = \{f \in \mathcal{F}: \text{supp} (\hat{f}) \subseteq \theta \}.
            \end{equation*}
            Note that each $f \in \mathcal{F}_{\theta}$ is such that $\norm{f_{\theta}}_{L^{\infty}(\R^2)} \leq 1$, and every $f \in \mathcal{F}$ is such that $f= \sum_{\theta} f_{\theta}$, for some $f_{\theta} \in \mathcal{F}_{\theta}$. A corollary of the contraction principle for Rademacher processes\footnote{Alternatively, one could compare this process to a Gaussian process and use Slepian's lemma.} (see \cite[Exercise 7.4]{vanHandel2016}, or Appendix B) is that
            \begin{equation}
                \mathbb{E}_{\epsilon} \sup_{f \in \mathcal{F}} \abs{\sum_{k \in S} \epsilon_k  \abs{D^{(\alpha)}f(y_{kl})}} \lesssim \mathbb{E}_{\epsilon}\sup_{f \in \mathcal{F}}\abs{\sum_{k \in S} \epsilon_k D^{(\alpha)}f(y_{kl})} .
            \end{equation}
            We may write\footnote{Notice that the fact $f \rightarrow D^{(\alpha)}f$ is a \textit{linear map} is crucial in the following step. See Section \ref{furtherremarkssection} for further discussion.}
            \begin{align}\label{thetareductionest}
                \mathbb{E}_{\epsilon}\underset{\substack{f \in \mathcal{F} \\ f= \sum f_\theta}}{\sup}\abs{\sum_{\theta} \sum_{k \in S} \epsilon_k D^{(\alpha)}f_\theta(y_{kl})} &\leq \sum_\theta \mathbb{E}_{\epsilon}\underset{\substack{f \in \mathcal{F} \\ f= \sum f_\theta}}{\sup}\abs{ \sum_{k \in S} \epsilon_k D^{(\alpha)}f_\theta(y_{kl})} \notag\\
                &\leq \sum_\theta \mathbb{E}_{\epsilon} \sup_{f \in \mathcal{F}_\theta} \abs{ \sum_{k \in S} \epsilon_k D^{(\alpha)}f(y_{kl})},
            \end{align}
        since $\abs{\{\theta\}} = N$, in light of (\ref{thetareductionest}), it suffices to establish the estimate
        \begin{equation}\label{finalboss}
            \mathbb{E}_{\epsilon} \sup_{f \in \mathcal{F}_{\theta}}\abs{\sum_{k \in S} \epsilon_k  D^{(\alpha)}f(y_{kl})} \lesssim_{\epsilon} N^{1/2+\epsilon} \sqrt{\abs{S}}, 
        \end{equation}
        for all $\theta$. 

        In order to establish (\ref{finalboss}), we require two ingredients. First, we require a description of elements in the class $\mathcal{F_{\theta}}$. This description will allow us to control the complexity of the set $\mathcal{F}_{\theta}$. Such a description is referred to as a wave packet decomposition. The basic idea is that if $\text{supp}(\hat{f_{\theta}}) \subseteq \theta$, and $\norm{f_{\theta}}_{L^{\infty}} \leq 1$, then when $f_{\theta}$ is restricted to $B_{N^2}$ we can write
        \begin{equation}\label{wavepacketheuristic}
            f_{\theta} \approx \sum_{k=1}^{N} a_k \mathbbm{1}_{\theta_{k}^{*}},
        \end{equation}
        where $\abs{a_k} \leq 1$ for all $k$ and the $\theta_k^{*}$'s are $N \times N^2$ rectangles oriented in the dual direction to $\theta$ that cover $B_{N^2}$. The wave packet decomposition plays an important role in several problems in geometric Fourier analysis, see for example, \cite{Guth2016} or \cite{CarberyIliopoulouWang2024} for precise statements and applications. Of course, (\ref{wavepacketheuristic}) is not literally true, and so we need the following precise version for our purposes, which is proved in the appendix.
        \begin{proposition}[Abstract Wave Packet Decomposition]\label{wavepacket}
                Suppose $f \in \mathcal{S}(\R^2)$ is such that $ \text{supp}(\hat{f}) \subseteq \theta$ and $\norm{f}_{L^{\infty}(\R^2)} \lesssim 1$. For every $\delta>0$, we can write
                \begin{equation*}
                    f(x) = \sum_{k \in \Lambda_{\delta}} a_k \phi_k(x) + R_f(x),
                \end{equation*}
                where the $\phi_k$'s are independent of the choice of $f$ and
                \begin{enumerate}
                    \item $\abs{a_k} \lesssim 1$ for all $k \in \Lambda_{\delta}$
                    \item $\abs{\Lambda_{\delta}} \lesssim N^{1+\delta}$
                    \item $\sup_{x \in \R^2} \sum \abs{\phi_k(x)} \lesssim 1$
                    \item $\norm{R_f(x)}_{L^{\infty}\left(B_{N^2}\right)} \lesssim_{\delta} N^{-M}$ for every $M>0$.
                \end{enumerate}
            \end{proposition}
            \begin{remark}
                We call this an `abstract' wave packet decomposition because for our purposes we do not need any geometric information about the $\phi_k$'s. Our proof does not require any ideas coming from incidence geometry or the geometry of wave packets. For more commentary on this matter, refer to the Appendix A.
            \end{remark}
            The second ingredient we require is a corollary of Dudley's inequality (see Appendix B for more background).
            \begin{proposition}\label{dudleyconsequence}
                Fix $S \subseteq [CN^4]$, and define the metric on the set $\mathcal{F}_{\theta}$ given as follows
                \begin{equation*}
                    d(f,g) = \left(\sum_{k \in S} \abs{Df^{(\alpha)}f(y_{kl}) - Dg^{(\alpha)}(y_{kl})}^2 \right)^{1/2}.
                \end{equation*}
                Let $\mathcal{N}(\mathcal{F}_\theta, d, u)$, be equal to the minimum number of balls of radius $u$ in the metric $d(\cdot, \cdot)$ needed to cover $\mathcal{F_{\theta}}$. Then,
                \begin{equation}\label{dudley}
                    \mathbb{E}_{\epsilon} \sup_{f \in \mathcal{F}_{\theta}}\abs{\sum_{k \in S} \epsilon_k  \abs{D^{(\alpha)}f(y_{kl})}} \lesssim \int_{0}^{\text{diam} (\mathcal{F}_\theta, d)} \sqrt{\log \mathcal{N} (\mathcal{F}_\theta, d, u)} du. 
                \end{equation}
            \end{proposition}
            For more background on this estimate refer to the appendix.
            \subsection{Entropy Estimates}
            In view of (\ref{dudley}), we need to obtain bounds on the quantity $\log \mathcal{N}(\mathcal{F}_\theta, d, u)$. We have the following proposition:
            \begin{proposition}\label{entropyestimate}
                For every $A>4$, there exist $c, C >0$ (possibly depending on $\alpha$) such that following estimate holds:
                \begin{align*}
                    \log \mathcal{N}(\mathcal{F}_{\theta}, d, u\sqrt{\abs{S}}) \lesssim_{\epsilon}
                    \begin{cases}
                        N^{1+\epsilon} \log (\frac{C}{u}) \quad \text{if  } u \geq cN^{-A}\\
                        \abs{S} \log (\frac{C}{u}) \quad \text{if  } 0<u < cN^{-A}
                    \end{cases}
                \end{align*}
                 \end{proposition}
                We first show that this proposition implies (\ref{finalboss}). 
                \begin{proof}[Proof of (\ref{finalboss}) assuming Proposition \ref{entropyestimate}]
                    First note by Bernstein's inequality, for any $f \in \mathcal{F}_{\theta}$,
                    \begin{equation}\label{basicbernstein}
                        \norm{D^{(\alpha)}f}_{L^{\infty}} \lesssim \norm{f}_{L^{\infty}} \lesssim 1.
                    \end{equation}
                    Therefore, 
                    \begin{equation}
                        d(f,g) = \left(\sum_{k \in S} \abs{Df^{(\alpha)}f(y_{kl}) - Dg^{(\alpha)}(y_{kl})}^2 \right)^{1/2} \lesssim \sqrt{\abs{S}},
                    \end{equation}
                    hence $\text{diam}(\mathcal{F_\theta}, d) \lesssim \sqrt{\abs{S}}$. 
                    Now
                    \begin{align*}
                        \mathbb{E}_{\epsilon} \sup_{f \in \mathcal{F}_{\theta}}\abs{\sum_{k \in S} \epsilon_k  \abs{D^{(\alpha)}f(y_{kl})}} &\lesssim \int_{0}^{C_1\sqrt{\abs{S}}} \sqrt{\log \mathcal{N} (\mathcal{F}_\theta, d, u)} du\\
                        &= \sqrt{\abs{S}}\int_{0}^{C_1} \sqrt{\log \mathcal{N} \left(\mathcal{F}_\theta, d, u \sqrt{\abs{S}} \right)} du.
                    \end{align*}
                    We estimate this integral in two parts using Proposition \ref{entropyestimate}. We have
                    \begin{align*}
                        &\int_{0}^{C_1} \sqrt{\log \mathcal{N} \left(\mathcal{F}_\theta, d, u \sqrt{\abs{S}} \right)} du\\
                        = &\int_{0}^{cN^{-A}} \sqrt{\log \mathcal{N} \left(\mathcal{F}_\theta, d, u \sqrt{\abs{S}} \right)} du + \int_{cN^{-A}}^{C} \sqrt{\log \mathcal{N} \left(\mathcal{F}_\theta, d, u \sqrt{\abs{S}} \right)} du  
                    \end{align*}
                    To estimate the first integral, note that $\abs{S} \lesssim N^4$, and therefore the first integral is bounded by
                    \begin{equation*}
                        \int_{0}^{cN^{-A}}\sqrt{\abs{S}} \cdot \sqrt{\log \left(\frac{C}{u} \right)} du \lesssim \sqrt{\abs{S}} N^{-A} \lesssim 1
                    \end{equation*}
                    since $\sqrt{\log(1/u)} \lesssim u^{-c}$ for all $0<c<1$ and $0<u<1$. Choosing $c$ sufficiently small, and $A> 4$ yields the claim. The second integral is bounded by
                    \begin{equation*}
                        \int_{0}^{C_1} N^{1/2+\epsilon} \sqrt{\log (C/u)} du \lesssim N^{1/2+ \epsilon}.
                    \end{equation*}
                        by the integrability of $\log(1/u)$ at $0$.
                        
                        Putting everything together we see that 
                        \begin{equation*}
                            \int_{0}^{C_1\sqrt{\abs{S}}} \sqrt{\log \mathcal{N} (\mathcal{F}_\theta, d, u)} du \lesssim_{\epsilon} \sqrt{\abs{S}} N^{1/2+ \epsilon}
                        \end{equation*}
                        which completes the proof.
                        \end{proof}

                        We are now ready to establish the entropy estimates in Proposition \ref{entropyestimate}. We begin with an elementary lemma on how metric entropy estimates are transfered under the image of Lipschitz map.
                        \begin{lemma}[Lipschitz Distortion of Metric Entropy]\label{lipschitzentropy}
                            Let $(X, d_{X})$ and $(Y, d_{Y})$ be two metric spaces and let $\Phi: X \rightarrow Y$ be a map such that $d_{X}(x_1, x_2) \leq L d_{Y}(\Phi(x_1),\Phi(x_2))$ for all $x_1, x_2 \in X$, then
                        \begin{equation*}
                              \log\mathcal{N}(X, d_{X}, Lu) \leq \log\mathcal{N}(\Phi(X), d_{Y}, u)
                        \end{equation*}
                        \end{lemma}
                        \begin{proof}
                            Let $\{\Phi(x_1), \Phi(x_2),..., \Phi(x_r)\} \in Y$ be a set of distinct elements such that for any $\Phi(x) \in \Phi(X)$, there is some $\Phi(x_i)$ for which $d_{Y}(\Phi(x), \Phi(x_i)) \leq u$. Then the set $\{x_1, x_2,...,x_r \}$ is such that if $x \in X$ is arbitrary, then there is some $x_i$ such that $d_{X}(x_i, x) \leq Ld_Y(\Phi(x_i), \Phi(x)) \leq Lu$. It follows that
                        \begin{equation*}
                            \mathcal{N}(X, d_{X}, Lu) \leq  \mathcal{N}(X, d_{Y}, u),
                        \end{equation*}
                        taking logarithms completes the proof.
                        \end{proof}
                        \begin{proof}[Proof of Proposition \ref{entropyestimate}]
                            Any $f \in \mathcal{F}_\theta$, and any $\delta>0$, by Proposition \ref{wavepacket} we may write 
                            \begin{equation}\label{repDf}
                                D^{(\alpha)}f(x) = \sum_{r=1}^{J} a_r \phi_r(x) + R_f(x)
                            \end{equation}
                            where the $\phi_r$'s and $R(x)$ satisfy properties 1-4 in Proposition \ref{wavepacket} and $J \lesssim N^{1+\delta}$.Proposition \ref{wavepacket} applies because $D^{(\alpha)}f$ is Fourier supported in $\theta$ and $\norm{D^{(\alpha)}f}_{L^{\infty}} \lesssim 1$, by (\ref{basicbernstein}). Consider the map $\Phi: \mathcal{F}_{\theta} \rightarrow \mathbb{C}^{J}$ given by
                            \begin{equation*}
                                \Phi(f) = (a_1,a_2,..,a_J),
                            \end{equation*}
                            where the tuple $(a_1,a_2,...,a_J)$ comes from (\ref{repDf}). Note that a representation of the form (\ref{repDf}) may not be unique, but we choose the $a_j$'s that come out of the construction of Proposition \ref{wavepacket} and its proof. Since each $\abs{a_k} \lesssim 1$, we note that
                            \begin{equation*}
                                \Phi(\mathcal{F_{\theta}}) \subseteq B_{\infty}^{C} \coloneqq \left\{(c_i)_{i=1}^{J} \subseteq \C^J: \abs{c_i} \leq C, \quad 1 \leq i \leq J\right\},
                            \end{equation*}
                            for some absolute constant $C$.
                            Let $f, g \in \mathcal{F}_\theta$, be arbitrary, and suppose that
                            \begin{equation*}
                                \Phi(f) = (a_1,a_2,...,a_J)  \quad \text{and} \quad \Phi(g)=(b_1,b_2,...,b_J).
                            \end{equation*}
                        We can write
                        \begin{equation}\label{lossy}
                            d(f,g) = \left(\sum_{k \in S}\abs{D^{(\alpha)}f(y_{kl}) - D^{(\alpha)}g(y_{kl})}^2\right)^{1/2} \leq \sqrt{\abs{S}}\max_{k \in S} \abs{D^{(\alpha)}f(y_{kl}) - D^{(\alpha)}g(y_{kl})}
                        \end{equation}
                        Using the representation (\ref{repDf}) (for $f$ and $g$) we have
                        \begin{align*}
                            \max_{k \in S} \abs{D^{(\alpha)}f(y_{kl}) - D^{(\alpha)}g(y_{kl})} &= \max_{k \in S} \abs{\sum_{r=1}^{J} (a_r-b_r) \phi_r(y_{kl}) + R_f(y_{kl}) - R_g(y_{kl})}\\
                            &\lesssim \max_{k \in S}  \abs{\sum_{r=1}^{J} (a_r-b_r) \phi_r(y_{kl})} + N^{-A + 1 + \delta}\\
                            &\lesssim \norm{a-b}_{\ell_{\infty}(\mathbb{C}^J)} \cdot \max_{k \in S}  \left(\sum_{k\in S} \abs{\phi_k(y_{kl})} \right) +N^{-A+1+ \delta}\\
                            &\lesssim \norm{a-b}_{\ell_\infty(\mathbb{C}^{J})} + N^{-A + 1 + \delta}
                        \end{align*}
                        where the second and third inequalities follows from properties 4 and 3 of Proposition \ref{wavepacket} respectively. Thus when $\norm{a-b}_{\ell^{\infty}} \gtrsim N^{-A}$ (where $A$ can be taken arbitrarily large), we have
                        \begin{equation*}
                            d(f,g) \lesssim \sqrt{\abs{S}} \norm{a-b}_{\ell_{\infty}(\mathbb{C}^J)}
                        \end{equation*}
                        Therefore by Lemma \ref{lipschitzentropy}, when $u \gtrsim N^{-A}$, 
                        \begin{equation}
                            \log(\mathcal{F}_{\theta}, d, u) \leq \log \mathcal{N}(B^C_\infty, \ell_\infty, u).
                        \end{equation}
                        It turns out that estimating the term $\log \mathcal{N}(B^C_\infty, \ell_\infty, u)$ is quite easy, it essentially is asking for how many copies of a small ball in the $\ell_\infty$ norm are are needed to cover a large ball in the the $\ell_\infty$ norm. Here is a sketch of this estimate. The $\ell_{\infty}$ ball of radius $C$ in $\R^{2J}$, can be covered perfectly by $\ell_{\infty}$ balls of radius $u$. Comparing the volumes of these two bodies, and computing their quotient we see that
                        \begin{equation}\label{linftyentropysimple}
                            \mathcal{N}(B^C_\infty, \ell_\infty, u) \lesssim \left(\frac{C}{u}\right)^{2J}.
                        \end{equation}
                        Finally, note that the complex $\ell_{\infty}$ ball in $\mathbb{C}^J$ is contained in the real $\ell_{\infty}$ ball in $\mathbb{R}^{2J}$, which yields the desired bound
                        \begin{equation*}
                            \log \mathcal{N}(B^C_\infty, \ell_\infty, u) \lesssim N^{1+\delta} \log \left(\frac{C}{u}\right),
                        \end{equation*}
                        and choosing $\delta$ sufficiently small proves first part of Proposition \ref{entropyestimate}. In order to prove the second part of the proposition, we forget about the wave packet structure of elements in $\mathcal{F}_\theta$ and proceed more crudely. For every $f \in \mathcal{F}_{\theta}$, consider the map $\Phi': \mathcal{F_{\theta}} \rightarrow \mathbb{C}^{\abs{S}}$
                        \begin{equation}
                            \Phi'(f) = (D^{(\alpha)}f(y_{kl}))_{k \in S}
                        \end{equation}
                        Since $\norm{D^{(\alpha)}f}_{L^\infty} \lesssim 1$, proceeding as before we see that
                        \begin{equation*}
                            d(f,g) \lesssim \sqrt{\abs{S}} \cdot \norm{\Phi'(f)-\Phi'(g)}_{\ell_\infty(\C^\abs{S})}, 
                        \end{equation*}
                        and 
                        \begin{equation}
                            \Phi'(\mathcal{F}_\theta) \subseteq B'^{C}_{\infty} \coloneqq \left\{(c_i)_{i=1}^{J} \subseteq \C^{\abs{S}}: \abs{c_i} \leq C, \quad 1 \leq i \leq \abs{S}\right\}.
                        \end{equation}
                        With the exact same argument used to produce the estimate (\ref{linftyentropysimple}), we see that
                        \begin{equation*}
                            \mathcal{N}(B_{\infty}^{'C}, \ell_\infty, u) \lesssim \left(\frac{C}{u}\right)^{\abs{S}},
                        \end{equation*}
                        so that 
                        \begin{equation*}
                            \log \mathcal{N}((B_{\infty}^{'C}, \ell_\infty, u) \lesssim \abs{S}\log \left(\frac{C}{u}\right),
                        \end{equation*}
                        which completes the proof of the proposition.
                        \end{proof}
\section{Proof of the Lower Bound}\label{lowerboundsection}
In this section, we will show that the estimate (\ref{mainpropeq}) is sharp by constructing a family of wave packets that saturate this bound.  For a fixed cap $\theta$, let $\mathbb{T}_\theta$ be the collection of $N \times N^2$ rectangles that are oriented in the dual direction to $\theta$ that cover $B_{N^2}$.
\begin{proposition}
    Fix $\theta$. There exist functions $\{\phi_{T}\}_{T \in \mathbb{T}_\theta}$ with each $\phi_{T} \in \mathcal{F}_\theta$  with following properties
    \begin{enumerate}
        \item $\text{supp}(\hat{\phi}_T) \subseteq \theta$
        \item $\abs{\phi_{T}} \gtrsim 1$ on $\frac{1}{10} T$ for each $T \in \mathbb{T}_\theta$.
        \item $\sup_{x \in B_{N^2}} \sum_{T \in \mathbb{T}_\theta} \abs{\phi_T (x)} \lesssim 1$.
    \end{enumerate} 
\end{proposition}
Such a construction is well known and we do not prove it here. We refer the interested reader to \cite[Notes 3]{Guth2017Decoupling}. 

Let $S$ be the (random) set
\begin{equation}
    S \coloneqq \{k: \delta_k = 1\},
\end{equation}
and following (\ref{randomsetballs}), let $V$ be the random set of balls given by
\begin{equation*}
    V \coloneqq \left\{ \bigcup_{k: \delta_k =1} V_{k} \right\} \subseteq B_{N^2}.
\end{equation*}
Suppose $V_k \subseteq V$ is one of the randomly chosen balls. We say $T \in \mathbb{T}_\theta$ 
 intersects $V_k$ \textit{non-trivially} if $\abs{V_k \cap T } \geq \frac{1}{2}\abs{V_k}$. We record the following basic geometric observations: For any fixed $V_k$, and each cap $\theta$, there is at least one $T \in \mathbb{T}_\theta$ that non-trivially intersects $V_k$. Moreover, for every $T \in \mathbb{T}$ that non-trivially intersects $V_k$, there is $V_k' \subseteq V_k$ such that $\abs{V_k'} \sim \frac{1}{100} \abs{V_k}$ and $\abs{V_k'} \subseteq T$ for \textit{all} such $T$. A consequence of this is that
 \begin{equation*}
     \abs{\int_{V_k'} \phi_T(x)dx} \gtrsim 1,
 \end{equation*}
for every $T$ that non-trivally intersects $V_k$.
Let $\mathbb{T} = \bigcup_\theta \mathbb{T}_\theta$. For each $T \in \mathbb{T}$, let $c_T$ be such that $\abs{c_T} \leq 1$ so that the function
\begin{equation}\label{wavepacketrepresentation}
    f = \sum_{T \in \mathbb{T}} c_T \phi_T
\end{equation}
is such that $f \in \mathcal{F}$.
 By the bounded overlap of the $V_k$'s we have
\begin{align*}
    \int_V \abs{f(x)} dx &\gtrsim \sum_{k \in S} \int_{V_k} \abs{f(x)}dx
    \geq \sum_{k \in S} \abs{\int_{V_k} f(x) dx}
\end{align*}
We choose the $c_k$'s to be independent random signs $(\epsilon_T)_{T \in \mathbb{T}}$ for a particular choice of the $c_T$'s and use Khintchine's inequality. More precisely, we have
\begin{equation*}
    \sup_{f \in \mathcal{F}} \int_{V} \abs{f(x)}dx \gtrsim \sum_{k \in S} \int_{V_k}\mathbb{E}_{\epsilon} \abs{\sum_{T \in \mathbb{T}} \epsilon_j \phi_T(x)} dx \sim \sum_{k \in S} \int_{V_k} \left(\sum_{T \in \mathbb{T}} \abs{\phi_T(x)}^2\right)^{1/2}dx.
\end{equation*}
For any fixed $k$ and fixed $\theta$, note that there is atleast one $T \in \mathbb{T}_\theta$ that non trivially intersects $V_k$. Since there are $N$ caps $\theta$, we see that
\begin{equation*}
    \int_{V_k} \left(\sum_{T \in \mathbb{T}} \abs{\phi_T(x)}^2\right)^{1/2}dx \gtrsim N^{1/2},
\end{equation*}
for every $k$. Thus
\begin{equation}
    \mathbb{E} \sup_{f \in \mathcal{F}}\int_{V}\abs{f(x)}dx \gtrsim N^{1/2} \mathbb{E} \abs{S} = N^{1/2}q
\end{equation}
which proves one part of the lower bound. To proceed with the second part of the lower bound, first note that
\begin{equation*}
    \int_V \abs{f(x)} dx \gtrsim \sum_{k \in S} \int_{V_k'} \abs{f(x)}dx
    \geq \sum_{k \in S} \abs{\int_{V_k'} f(x) dx}, 
\end{equation*}
letting $f$ be chosen as in (\ref{wavepacketrepresentation}), we obtain
\begin{equation}
    \sup_{f \in \mathcal{F}
    }\int_{V}\abs{f(x)}dx \gtrsim \sum_{k \in S}\abs{ \sum_{T \in \mathbb{T}} \int_{V_k'} c_T \phi_T(x) dx},
\end{equation}
where we are free to choose the $c_T$'s as long as $\abs{c_T} \leq 1$ for all $T$. Let $(\epsilon_k)_{k \in S}$ be an independent and identically distributed collection of Rademacher signs. The previous inequality gives
\begin{equation*}
    \sup_{f \in \mathcal{F}
    }\int_{V}\abs{f(x)}dx \gtrsim \mathbb{E}_{\epsilon} \sum_{k \in S}\abs{ \epsilon_k \sum_{T \in \mathbb{T}} \int_{V_k'} c_T \phi_T(x) dx} \geq \mathbb{E}_{\epsilon} \abs{\sum_{T \in \mathbb{T}}c_T \sum_{k \in S}\epsilon_k \int_{V_k'} \phi_T(x) dx}.
\end{equation*}
Since we are free to choose the $c_T$'s, for each $T \in \mathbb{T}$, we choose $c_T$ so that
\begin{equation*}
    c_T \left(\sum_{k \in S} \epsilon_k \int_{V_k'}\phi_T(x)dx \right) = \abs{\sum_{k \in S} \epsilon_k \int_{V_k'}\phi_T(x)dx}.
\end{equation*}
Therefore,
\begin{equation*}
    \sup_{f \in \mathcal{F}} \int_{V}\abs{f(x)}dx \gtrsim \sum_{T \in \mathbb{T}}\mathbb{E}_{\epsilon} \abs{\sum_{k \in S} \epsilon_k \int_{V_k'}\phi_T(x)dx} \sim \sum_{T \in \mathbb{T}}\left(\sum_{k \in S}\abs{\int_{V_k'}\phi_T(x)dx}^2\right)^{1/2}
\end{equation*}
by Khintchine's inequality. For any fixed $T \in \mathbb{T}$, note that $T$ contains $\sim N^3$ unit balls. The probability of any unit ball $V_k \subset T$ to be included in $V$ is $\sim \frac{q}{N^4}$. Thus for any fixed $T$ we observe that  
\begin{equation*}
    \sum_{k \in S}\abs{\int_{V_k'}\phi_T(x)dx}^2 \sim \text{Binomial}\left(N^3, \frac{q}{N^4}\right),
\end{equation*}
which is a random variable with mean $\sim \frac{q}{N}$. Finally,
\begin{equation*}
    \mathbb{E}_{S} \sup_{f \in \mathcal{F}} \int_{V}\abs{f(x)}dx \gtrsim \mathbb{E}_{S} \sum_{T \in \mathbb{T}}\left(\sum_{k \in S}\abs{\int_{V_k'}\phi_T(x)dx}^2\right)^{1/2} \sim N^2 \cdot \left(\frac{q}{N}\right)^{1/2} = N^{3/2}\sqrt{q},
\end{equation*}
which completes the proof of the second part of the lower bound \footnote{Note that in the final step we use the following fact: If $X \sim \text{Binomial}(N^3, \frac{q}{N^4})$, then $\mathbb{E}\sqrt{X} \sim \sqrt{\mathbb{E}(X)}$. There are several ways to prove this. One way to do this is via Taylor approximation, see \url{https://mathoverflow.net/questions/121411/expectation-of-square-root-of-binomial-r-v} , which shows that this claim is true when $np \gtrsim 1$, which is true in our case since $q \geq N$. }.

This lower bound proof bears similarity with the `wave packet tuning' strategy presented in \cite{guth2022enemy}. I am sure this is known to experts, but I have not seen this precisely written down anywhere in the literature.
\section{Further Remarks}\label{furtherremarkssection}
\subsection{Suggestions For Reading This Paper}
This paper makes several technical reductions that make it possible for the Fourier analytic setup for the problem to work with probabilistic machinery. On the first reading of this paper, the reader is suggested to make the following `cheats' to understand the basic strategy of the proof.
\begin{enumerate}
    \item Assume that every $f \in \mathcal{F}$ is constant on any unit ball $V_k \subseteq B_{N^2}$. If this is assumed, the technical reduction of Section \ref{reductionsection} is not needed. One can replace the term $D^{(\alpha)}f(y_{kl})$ with $f(y_k)$.
    \item Assume that any $f \in \mathcal{F}_\theta$ can be written as
    \begin{equation*} 
        f = \sum_{T \in \mathbb{T}_\theta}c_T \mathbbm{1}_T
    \end{equation*}
    where $\mathbb{T}_{\theta}$ is the set of $N \times N^2$ rectangles that cover $B_{N^2}$. This avoids the need for technical estimates used to prove Proposition \ref{wavepacket}.
    \item As a further simplication, it may be helpful to consider the special case where we restrict to the class of $f \in \mathcal{F}$, where $f$ is the sum of \textit{only} one wave packet in each direction. That is, consider the simplified model where each $f \in \mathcal{F}$ can be represented as 
    \begin{equation*}
        f = \sum_{\theta} c_\theta \mathbbm{1}_{T_\theta}
    \end{equation*}
    for some $T_\theta \in \mathbb{T}_\theta$. In this case one obtains the bound
    \begin{equation}\label{squarerootcancellation}
        \sup_{f \in \mathcal{F}} \norm{f}_{L^1(V)} \sim  N^{1/2}q,
    \end{equation}
    that is, the first term in Proposition \ref{mainprop} dominates. 
\end{enumerate}
On a personal note, when I started working on this problem, I conjectured that (\ref{squarerootcancellation}) was the truth. I was not able to prove this, but eventually managed to do so in the simplified model of case (3) mentioned above. The argument in this proof (in particular the inequality (\ref{lossy})) feels quite lossy, and I was not convinced that this argument was sharp until I discovered the second part of the lower bound in Section \ref{lowerboundsection}.

Another important detail is why the reduction of Section \ref{reductionsection} is necessary. As far as I can tell, the only place where it is necessary in the sequence of inequalities culminating in (\ref{thetareductionest}). In order for this sequence of inequalities to work one crucially needs that the map
\begin{equation*}
    f \rightarrow D^{(\alpha)}f
\end{equation*}
is a \textit{linear} map. If we do not have the reduction of Section \ref{reductionsection}, then manipulating the map
\begin{equation*}
    f \rightarrow \int_{V_k}\abs{f(x)} dx 
\end{equation*}
is awkward, because it is non-linear and one cannot deduce the equivalent of (\ref{thetareductionest}) in this setting.
\newpage        
\section{Appendix A: Some Results in Fourier Analysis}
    \subsection{Wave Packet Decomposition}
            In this section we rigorously prove Propsition \ref{wavepacket}. The wave packet decomposition has been proved and used in various forms (see \cite{Guth2016}, \cite{CarberyIliopoulouWang2024}, \cite{GuthHickmanIliopoulou2019} for illustrative examples). In order to make our proof work we require $L^\infty$ control of the remainder term and also need a pointwise bound for the sum of the wave packets. The wave packet decompositions proved in the existing literature do not  explicitly prove these properties, so we prove this precise version for our purposes here. Moreover, in our case we do not require any geometric properties of the wave packets. Most proofs involving the wave packet decomposition need to show (precisely) that the wave packets are large on  a $N \times N^2$ rectangle and rapidly decaying outside. Our proof does not need this feature of the wave packet decomposition. 
            \begin{proposition}[Abstract Wave Packet Decomposition]
                Suppose $f \in \mathcal{S}(\R^2)$ is such that $ \text{supp}(\hat{f}) \subseteq \theta$ and $\norm{f}_{L^{\infty}(\R^2)} \lesssim 1$. For every $\delta>0$, we can write
                \begin{equation*}
                    f(x) = \sum_{k \in \Lambda_{\delta}} a_k \phi_k(x) + R_f(x),
                \end{equation*}
                where the $\phi_k$'s are independent of the choice $f$ and
                \begin{enumerate}
                    \item $\abs{a_k} \lesssim 1$ for all $k \in \Lambda_{\delta}$
                    \item $\abs{\Lambda_{\delta}} \lesssim N^{1+\delta}$
                    \item $\sup_{x \in \R^2} \sum_{k \in \Lambda_\delta} \abs{\phi_k(x)} \lesssim 1$
                    \item $\norm{R(x)}_{L^{\infty}\left(B_{N^2}\right)} \lesssim_{\delta} N^{-M}$ for every $M>0$.
                \end{enumerate}
            \end{proposition}
            \begin{proof}
                By rotating and translating $f$ in Fourier space (these operations do not affect the desired properties of our abstract decomposition, as can be checked), we may without loss of generality assume that $\theta = \left[0,\frac{1}{N}\right] \times \left[0,\frac{1}{N^2}\right]$. Define
                \begin{equation}
                    h(\xi_1, \xi_2) = \hat{f}\left(\frac{\xi_1}{N}, \frac{\xi_2}{N^2}\right).
                \end{equation}
                By Fourier inversion we can write
                \begin{equation}\label{fourierinversiontheta}
                    f(x) = \int\int_{\theta} \hat{f}(\xi_1,\xi_2) e^{2\pi i (\xi_1 x_1 + \xi_2 \cdot x_2)}d \xi.
                \end{equation}
         
            Applying the change of variables $\xi_1 \rightarrow \frac{\xi_1}{N}$ and $\xi_2 \rightarrow \frac{\xi_2}{N^2}$ in (\ref{fourierinversiontheta}), the region $\theta$ is transformed to the region $[0,1]^2$ and the Jacobian factor is $\frac{1}{N^3}$. Therefore we have
            \begin{equation}\label{fourierinversion}
                f(x) = \frac{1}{N^3}\underset{[0,1]^2}{\int\int} \hat{f}\left(\frac{\xi_1}{N}, \frac{\xi_2}{N^2} \right) e^{2\pi i \left(\frac{\xi_1 x_1}{N}  + \frac{\xi_2 x_2}{N^2} \right)} d\xi = \frac{1}{N^3}\underset{[0,1]^2}{\int\int} h\left(\xi_1, \xi_2 \right) e^{2\pi i \left(\frac{\xi_1 x_1}{N}  + \frac{\xi_2 x_2}{N^2} \right)} d\xi.
            \end{equation}
            Since $h$ is smooth and $\text{supp} (h) \subseteq [-\pi,\pi]^2$, we can extend $h$ to a $2\pi$ periodic function on $\R^2$. Expanding $h$ as a Fourier series on $[-\pi, \pi]$, we write
            \begin{equation}\label{fourierseries}
                h(\xi_1,\xi_2) = \sum_{n,m \in \Z} c_{n,m} e^{i(\xi_1 n + \xi_2 m)},
            \end{equation}
            where the above sum converges uniformly on $[0,1]^2$ (since $f \in \mathcal{S}(\R^2)$) and the Fourier coefficients are given by
            \begin{equation*}
                c_{n,m} = \frac{1}{2\pi}\underset{[0,1]^2}{\int \int} h(\xi_1, \xi_2)e^{-i(n\xi_1 + m\xi_2)}d\xi.
            \end{equation*}
            Recalling that $h(\xi_1, \xi_2) = f(\frac{\xi_1}{N}, \frac{\xi_2}{N^2})$, we use the change of  variables $\xi_1 \rightarrow \xi_1N$ and $\xi_2 \rightarrow \xi_2 N^2$ to yield 
            \begin{equation*}
                c_{n,m} = \frac{N^3}{2 \pi} \int \int_{\theta} \hat{f}(\xi_1, \xi_2)e^{-i (\xi_1 n N + \xi_2 mN^2)} d \xi.
            \end{equation*}
            Applying the Fourier inversion formula we see that
            \begin{equation*}
                c_{n,m} = \frac{N^3}{2 \pi } f\left(-\frac{nN}{2 \pi}, -\frac{mN^2}{2\pi} \right).
            \end{equation*}
            the condition that $\norm{f}_{L^{\infty}} \lesssim 1$ shows that $\abs{c_{n,m}} \lesssim N^3$.
            Now let $\eta(\xi_1, \xi_2) \in C_{c}^{\infty}(\R^2)$ be such that $\eta \equiv 1$ on $[0,1]^2$ and $\eta \equiv 0$ on $\R^2 \setminus [-2,2]^2$. Since $\text{supp}(h) \subseteq [0,1]^2$, clearly we have that
            \begin{equation*}
                h(\xi_1,\xi_2)= \eta(\xi_1, \xi_2) h(\xi_1, \xi_2), \forall (\xi_1, \xi_2) \in \R^2.
            \end{equation*}
            Combining (\ref{fourierinversion}) and (\ref{fourierseries}), we have that
            \begin{align*}
                f(x) &=\frac{1}{N^3}\underset{[0,1]^2}{\int\int} h\left(\xi_1, \xi_2 \right) e^{2\pi i \left(\frac{\xi_1 x_1}{N}  + \frac{\xi_2 x_2}{N^2} \right)} d\xi\\ 
                &= \frac{1}{N^3}\int\int_{\R^2} h\left(\xi_1, \xi_2 \right) \eta(\xi_1, \xi_2) e^{2\pi i \left(\frac{\xi_1 x_1}{N}  + \frac{\xi_2 x_2}{N^2} \right)} d\xi \\
                &= \sum_{n,m} \frac{c_{n,m}}{N^3}\int \int_{\R^2} \eta(\xi_1, \xi_2) e^{i (n\xi_1 + m\xi_2)}\cdot e^{2\pi i \left(\frac{\xi_1 x_1}{N}  + \frac{\xi_2 x_2}{N^2} \right)} d \xi\\
                &= \sum_{n, m} \frac{c_{n,m}}{N^3} \hat{\eta}\left(-\frac{x_1}{N} - \frac{n}{2\pi}, -\frac{x_2}{N^2} - \frac{m}{2\pi}. \right)
            \end{align*}
            We set
            \begin{equation*}
                a_{n,m} = \frac{c_{n,m}}{N^3}
            \end{equation*}
            and 
            \begin{equation*}
                \phi_{n,m}(x_1,x_2) = \hat{\eta}\left(-\frac{x_1}{N} - \frac{n}{2\pi}, -\frac{x_2}{N^2} - \frac{m}{2\pi}. \right)
            \end{equation*}
            It is clear that the $\phi_{n,m}$'s are independent of $f$ and $\abs{a_{n,m}} \lesssim 1$ for all $n,m$. Since $\hat{\eta}$ is the Fourier transform a bump function, geometrically, the function $\phi_{m,n}$ can be thought of as being largest on a $N \times N^2$ rectangle centered at the point $(\frac{nN}{2\pi}, \frac{mN^2}{2\pi})$ and rapidly decaying outside of this rectangle. We use this basic heuristic to set up our abstract wave packet decomposition. 

            \textbf{Selecting the wave packets and the remainder: } Fix $\delta > 0$, and define
            \begin{equation*}
                \Lambda_{\delta} \coloneqq\{(m,n) \in \Z^2: \abs{n} \leq 10N^{\delta}, \abs{m} \leq 10N^{1+\delta}\}
            \end{equation*}
            so that we have $\abs{\Lambda_{\delta}} \lesssim N^{1+ 2\delta}$. We set
            \begin{equation*}
                R_f(x_1,x_2) \coloneqq \sum_{(n,m) \not\in \Lambda_{\delta}} a_{n,m}\phi_{n,m}(x_1,x_2).
            \end{equation*}
            We can now write
            \begin{equation}
                f(x) = \sum_{(n,m) \in \Lambda_{\delta}} a_{n,m} \phi_{n,m}(x) + R(x),
            \end{equation}
            and properties 1 and 2 of the proposition have been established. It now remains to prove properties 3 and 4.
            
            The Schwartz decay of $\hat{\eta}$ tells us that for all large $M>0$,
            \begin{equation}\label{schwartzdecay}
                \abs{\phi_{n,m}(x_1,x_2)} \lesssim \left(1 + \abs{ \frac{n}{2\pi} + \frac{x_1}{N}} + \abs{\frac{m}{2\pi} + \frac{x_2}{N^2}} \right)^{-M}.
            \end{equation}
            We can now show that
            \begin{equation*}
                \sup_{x \in \R^2} \sum_{(n,m) \in \Lambda_{\delta}} \abs{\phi_{n,m}(x)} \lesssim 1,
            \end{equation*}
            which is property 3 of our wave packet decomposition. First we note that given any $(x_1,x_2) \in \R^2$, define
            \begin{equation*}
                E_j = \left\{(n,m) \in \Z^2 : j \leq \text{dist}_{\ell_1} \left( \left(\frac{-n}{2\pi},\frac{-m}{2 \pi}\right), \left(\frac{x_1}{N},\frac{x_2}{N^2} \right)\right) \leq j+1 \right\}, 
            \end{equation*}
            a basic counting argument gives
            \begin{equation}\label{basiccounting}
                \abs{E_j} \lesssim j.
            \end{equation}
            Now by (\ref{schwartzdecay}) and (\ref{basiccounting}),
            \begin{align*}
                \sum_{(n,m) \in \Lambda_{\delta}} \abs{\phi_{n,m}(x_1,x_2)} &\lesssim \sum_{(n,m) \in \Lambda_{\delta}} \left(1 + \abs{ \frac{n}{2\pi} + \frac{x_1}{N}} + \abs{\frac{m}{2\pi} + \frac{x_2}{N^2}} \right)^{-M}\\
                &\lesssim \sum_{j}\sum_{\substack{(n,m) \in \Lambda_{\delta} \\ (n,m) \in E_j}} \left(1 + \abs{ \frac{n}{2\pi} + \frac{x_1}{N}} + \abs{\frac{m}{2\pi} + \frac{x_2}{N^2}} \right)^{-M}\\
                &\lesssim_{M} \sum j (1+j)^{-M}\\
                &\lesssim 1,
            \end{align*}
            where the last line follows from choosing $M>3$.
            
            We now show that $\norm{R_f(x)}_{L^{\infty}(B_{N^2})} \lesssim_{\delta} N^{-A}$, for all $A>0$, which is the final property that needs to be proved. If $(n,m) \not\in \Lambda_{\delta}$ then either
            \begin{equation*}
                \abs{n} > 10N^{\delta} \quad \text{or} \quad \abs{m} > 10N^{1+\delta}.
            \end{equation*}
            In any case if $(x_1,x_2) \in B_{N^2}$ ,then we must necessarily have that
            \begin{equation}\label{Ndeltaupperbound}
                \abs{\frac{n}{2\pi} + \frac{x_1}{N}} + \abs{\frac{m}{2\pi} + \frac{x_2}{N^2}} \gtrsim \text{dist}_{\ell_1}\left(\left(-\frac{n}{2\pi},\frac{-m}{2\pi} \right), \left(\frac{x_1}{N},\frac{x_2}{N^2} \right)\right) \gtrsim N^{\delta},
            \end{equation}
            since $\abs{x_1}, \abs{x_2} \leq N^2$. 
            Finally, we can estimate the remainder: If $(x_1,x_2) \in B_{N^2}$ we have by (\ref{schwartzdecay}), (\ref{basiccounting}), and (\ref{Ndeltaupperbound}) that 
            \begin{align*}
                \abs{R_f(x)} &= \abs{\sum_{ (n,m) \not\in \Lambda_{\delta}} a_{n,m} \phi_{n,m}(x_1,x_2)}\\
                &\lesssim \sum_{(n,m) \not\in \Lambda_{\delta}} \abs{\phi_{n,m}(x_1,x_2)} \\
                &\lesssim \sum_{(n,m) \not\in \Lambda_{\delta}} \left(1+ \abs{\frac{n}{2\pi} + \frac{x_1}{N^2}} + \abs{\frac{m}{2\pi} + \frac{x_2}{N^2}}\right)^{-M}\\
                &\lesssim \sum_{j} \sum_{\substack{(n,m) \in \Lambda_{\delta} \\ (n,m) \in E_j}} j \cdot (1+j)^{-M}\\
                &\lesssim \sum_{j \gtrsim N^{\delta}}j(1+j)^{-M}\\
                &\lesssim_{\delta} N^{-A}
            \end{align*}
            where the last line holds if $M$ (depending on $\delta$) is chosen large enough. This shows that $\norm{R(x)}_{L^{\infty}(B_{N^2})} \lesssim_{\delta} N^{-A}$ as desired.
            \end{proof}
    We now provide a proof of the estimate (\ref{localconstl1}), as a consequence for the following two lemmas. 
    {\begin{proposition}
    We have the uniform bound
    \begin{equation*}
        \sup_{f \in \mathcal{F}} \sum_{k=1}^{CN^4} \abs{D^{(\alpha)}f(y_{kl})} \lesssim_{\epsilon} N^{9/2}
    \end{equation*}
\end{proposition}}
 For each $f \in \mathcal{F}$, we write $D^{(\alpha)}f(y_{kl}) = \int f(x) \eta^{(\alpha)}(x-y_{kl})dx$. As a consequence, for every $f \in \mathcal{F}$ we have that 
 \begin{equation*}
     \sum_{k=1}^{CN^4} \abs{D^{(\alpha)}f(y_{kl})} \leq \int_{\R^2} \abs{f(x)} \cdot \left( \sum_{k=1}^{CN^4}\abs{\eta^{(\alpha)}(x-y_{kl})}dx \right).
 \end{equation*}
 We now estimate this integral based on the following two lemmas:
    \begin{lemma}
     The following point-wise bound holds: For each $x \in B_{2N^2}$, we have
     \begin{equation*}
         \sum_{k=1}^{CN^4} \abs{\eta^{(\alpha)}(x-y_{kl})} \lesssim 1. 
     \end{equation*}
 \end{lemma}
 \begin{proof}
          The Schwartz Decay of $\eta^{(\alpha)}$ tells us that
     \begin{equation*}
         \abs{\eta^{(\alpha)}(x-y_{kl})} \leq C_{N} \text{dist} (x, y_{kl})^{-N}
     \end{equation*}
     for some $N$ to be chosen later. Define for integers $j \geq 1$,
     \begin{equation*}
         C_{j}= \{ y_{kl}: 2^{j-1} < \text{dist}(y_{kl}, x) \leq 2^{j} \},
     \end{equation*}
     and 
     \begin{equation*}
         C_0 = \{ y_{kl}: \text{dist}(y_{kl}, x) \leq 2^{0}\}.
     \end{equation*}
Since the set $\{y_{kl}\}$ is a 1 - separated set, we have that $\abs{C_j} \lesssim 2^{j}$, then for any $N \geq 3$,
\begin{align*}
    \sum_{k=1}^{CN^4} \abs{\eta^{(\alpha)}(x-y_{kl})} &\lesssim \sum_{j=0}^{\infty} C_N \text{dist}(x,y_{kl})^{-N} \abs{C_j} \\
    &\lesssim \sum_{j=0}^{\infty} 2^{j}2^{-N(j-1)}\\
    &\lesssim 1.
\end{align*}
 \end{proof}
\begin{lemma}
    We have
    \begin{equation*}
        \sup_{f \in \mathcal{F}}\norm{f}_{L^1(B_N^2)} \lesssim_{\epsilon} N^{9/2+ \epsilon}.
    \end{equation*}
    \begin{proof}
        For every $f \in \mathcal{F}$, using the superlevel set bound (\ref{superlevel}), we have that
        \begin{align*}
            \norm{f}_{L^1(B_{N^2})} &\sim \int_{0}^{N^2}\abs{\{x \in B_{N^2}: \abs{f(x)}> 
            \lambda\}} d\lambda \\
            & \sim  \int_{0}^{N^{1/2}}\abs{\{x \in B_{N^2}: \abs{f(x)}> 
            \lambda\}} d\lambda + \int_{N^{1/2}}^{N^{2}}\abs{\{x \in B_{N^2}: \abs{f(x)}> 
            \lambda\}} d\lambda \\
            &\lesssim_{\epsilon} N^{1/2}(N^4) + N^{\epsilon} \int_{N^{1/2}}^{N^2} \frac{N^7}{\lambda^6} d\lambda\\
            & \sim_{\epsilon} N^{9/2+ \epsilon}.
        \end{align*}
    \end{proof}
\end{lemma}
\newpage
\section{Appendix B: Probability Background}
We give some basic background on the theory behind estimating the suprema of Rademacher processes. These results are well known, and we provide them here for the reader who may not have much of a background in probability. For more background, we refer the reader to the excellent references \cite{vershynin2018high}, \cite{talagrand2014upper}, and \cite{vanHandel2016}.

Let $(\epsilon_{i})_{i=1}^{n}$ be an i.i.d. collection of Rademacher random variables, and suppose $T \subseteq \C^n$. We are interested in the quantity 
\begin{equation}\label{rademacherprocess}
    \mathbb{E} \sup_{t \in T} \abs{\sum_{i=1}^{n}\epsilon_i t_i}, \quad t = (t_1,t_2,...t_n).
\end{equation}
A way to find an upper bound for the quantity (\ref{rademacherprocess}) is to use Dudley's inequality. For any $s, t \in T$, define the canonical distance between $s$ and $t$ as follows:
\begin{equation*}
    d(s,t) \coloneqq \left(\sum_{i=1}^{n}\abs{s_i-t_i}^2 \right)^{1/2}.
\end{equation*}
\begin{proposition}[Dudley's inequality]
We have the following inequality
\begin{equation}\label{dudley}
    \mathbb{E}\sup_{t \in T}\abs{\sum_{i=1}^{n}\epsilon_i t_i} \lesssim \int_{0}^{\infty}  \sqrt{\log \mathcal{N}(T,d, \epsilon)} d\epsilon.
\end{equation}
\end{proposition}
Dudley's inequality is proved using the chaining method and we refer the reader to \cite{vershynin2018high} for a self contained proof.

Another tool we need is the contraction principle for Rademacher processes.
\begin{theorem}[Contraction Principle] For any $T \subseteq \C^n$ the following always holds
    \begin{equation*}
        \mathbb{E}\sup_{t \in T}\abs{\sum_{i=1}^{n}\epsilon_i \abs{t_i}} \lesssim \mathbb{E}\sup_{t \in T} \abs{\sum_{i=1}^{n}\epsilon_i t_i }
    \end{equation*}
\end{theorem}
The intuition behind this theorem is that if the $\abs{t_i}$'s are all positive then the left hand side of the above equation could have cancellations due to the oscillatory nature of the $\epsilon_k$'s. The right hand side gets rid of the positivity of the $\abs{t_i}$'s and therefore is a larger quantity. For a proof see \cite[Exercise 7.4]{vanHandel2016}.   
\begin{proof}[Proof of Proposition \ref{dudleyconsequence}] For the process
    \begin{equation*}
        \mathbb{E}\sup_{f \in \mathcal{F}_\theta} \abs{\sum_{k \in S} \epsilon_k D^{(\alpha)}f(y_{kl})}, 
    \end{equation*}
    the canonical metric is
    \begin{equation*}
        d(f,g) = \left(\sum_{k \in S} \abs{D^{(\alpha)}f(y_{kl}) - D^{(\alpha)}g(y_{kl})}\right)^{1/2}. 
    \end{equation*}
    An application of Dudley's inequality (\ref{dudley}) precisely gives the statement of Proposition \ref{dudleyconsequence}.
\end{proof}
\newpage
\bibliographystyle{plain}
\bibliography{references}

@book{vershynin2018high,
  author       = {Vershynin, Roman},
  title        = {High-Dimensional Probability: An Introduction with Applications in Data Science},
  series       = {Cambridge Series in Statistical and Probabilistic Mathematics},
  volume       = {47},
  year         = {2018},
  publisher    = {Cambridge University Press},
  address      = {Cambridge},
  isbn         = {978-1-108-41519-4},
  doi          = {10.1017/9781108231596},
  url          = {https://doi.org/10.1017/9781108231596}
}

@article{cairo2025counterexample,
  title={A Counterexample to the Mizohata-Takeuchi Conjecture},
  author={Cairo, Hannah Mira},
  journal={arXiv preprint arXiv:2502.06137},
  year={2025},
  note={Submitted on February 10, 2025},
  url={https://arxiv.org/abs/2502.06137}
}

@misc{guth2022enemy,
  author       = {Larry Guth},
  title        = {An Enemy Scenario in Restriction Theory},
  year         = 2022,
  howpublished = {\url{https://www.youtube.com/watch?v=x-DET83UjFg}},
  note         = {Accessed: 2025-05-21}
}

@book{talagrand2014upper,
  author    = {Michel Talagrand},
  title     = {Upper and Lower Bounds for Stochastic Processes: Modern Methods and Classical Problems},
  publisher = {Springer},
  year      = {2014},
  series    = {Ergebnisse der Mathematik und ihrer Grenzgebiete. 3. Folge / A Series of Modern Surveys in Mathematics},
  volume    = {60},
  isbn      = {978-3-642-54074-5},
  doi       = {10.1007/978-3-642-54075-2},
  url       = {https://link.springer.com/book/10.1007/978-3-642-54075-2}
}

@article{bourgain1989bounded,
  author    = {Jean Bourgain},
  title     = {Bounded orthogonal systems and the $\Lambda(p)$-set problem},
  journal   = {Acta Mathematica},
  volume    = {162},
  number    = {3--4},
  pages     = {227--245},
  year      = {1989},
  doi       = {10.1007/BF02392838},
  url       = {https://doi.org/10.1007/BF02392838}
}

@article{talagrand1994majorizing,
  author    = {Michel Talagrand},
  title     = {Constructions of Majorizing Measures, Bernoulli Processes and Cotype},
  journal   = {Geometric and Functional Analysis},
  volume    = {4},
  number    = {6},
  pages     = {660--717},
  year      = {1994},
  doi       = {10.1007/BF01896658},
  url       = {https://link.springer.com/article/10.1007/BF01896658}
}

@misc{mulherkar2025random,
  title         = {Random Constructions for Sharp Estimates of Mizohata--Takeuchi Type},
  author        = {Mulherkar, Siddharth},
  year          = {2025},
  month         = jun,
  eprint        = {2506.05624},
  archivePrefix = {arXiv},
  primaryClass  = {math.CA},
  doi           = {10.48550/arXiv.2506.05624},
  url           = {https://arxiv.org/abs/2506.05624}
}

@article{mockenhauptSchlag2009majorantRandomSets,
  title   = {On the {Hardy}--{Littlewood} majorant problem for random sets},
  author  = {Mockenhaupt, G. and Schlag, W.},
  journal = {Journal of Functional Analysis},
  volume  = {256},
  number  = {4},
  pages   = {1189--1237},
  year    = {2009},
  doi     = {10.1016/j.jfa.2008.06.005}
}

@article{BourgainDemeter2015Decoupling,
  title   = {The proof of the {$\ell^2$} Decoupling Conjecture},
  author  = {Bourgain, Jean and Demeter, Ciprian},
  journal = {Annals of Mathematics},
  series  = {Second Series},
  volume  = {182},
  number  = {1},
  year    = {2015},
  pages   = {351--389},
  doi     = {10.4007/annals.2015.182.1.9}
}

@misc{CairoZhang2025PowerLossMizohataTakeuchi,
  title         = {Power loss for the Mizohata--Takeuchi conjecture on {$C^k$} convex hypersurfaces},
  author        = {Cairo, Hannah and Zhang, Ruixiang},
  year          = {2025},
  month         = dec,
  eprint        = {2512.08064},
  archivePrefix = {arXiv},
  primaryClass  = {math.CA},
  doi           = {10.48550/arXiv.2512.08064},
  url           = {https://arxiv.org/abs/2512.08064}
}

@article{CandesPlan2010RISPless,
  title   = {A probabilistic and RIPless theory of compressed sensing},
  author  = {Cand{\`e}s, Emmanuel J. and Plan, Yaniv},
  year    = {2010},
  eprint  = {1011.3854},
  archivePrefix = {arXiv},
  primaryClass  = {math.OC}
}

@article{IosevichMayeli2024LambdaSignalRecovery,
  title   = {Uncertainty Principles, Restriction, Bourgain's {$\Lambda_q$} Theorem, and Signal Recovery},
  author  = {Iosevich, Alex and Mayeli, Azita},
  journal = {Applied and Computational Harmonic Analysis},
  year    = {2024},
  volume  = {76},
  number  = {1},
  pages   = {101734},
  doi     = {10.1016/j.acha.2024.101734},
  url     = {https://doi.org/10.1016/j.acha.2024.101734}
}

@article{Tao2008UUP,
  title   = {The uniform uncertainty principle and compressed sensing},
  author  = {Tao, Terence},
  year    = {2008},
  url     = {https://terrytao.files.wordpress.com/2008/12/uup.pdf}
}

@misc{Guth2022DecouplingSurvey,
  author       = {Guth, Larry},
  title        = {Decoupling estimates in Fourier analysis},
  year         = {2022},
  note         = {Based on ICM 2022 survey lecture},
  eprint       = {2207.00652},
  archivePrefix= {arXiv},
  primaryClass = {math.CA},
  url          = {https://arxiv.org/abs/2207.00652}
}

@misc{GuthMaldagueWang2020ImprovedParabola,
  title         = {Improved decoupling for the parabola},
  author        = {Guth, Larry and Maldague, Dominique and Wang, Hong},
  year          = {2020},
  eprint        = {2009.07953},
  archivePrefix = {arXiv},
  primaryClass  = {math.CA},
  doi           = {10.48550/arXiv.2009.07953},
  url           = {https://arxiv.org/abs/2009.07953}
}

@misc{GuthYouTubeSharpExamples2023,
  author       = {Guth, Larry},
  title        = {Sharp examples for decoupling and related questions},
  howpublished = {YouTube video},
  year         = {2023},
  note         = {\url{https://www.youtube.com/watch?v=h_psKQV-pU4}},
  url          = {https://www.youtube.com/watch?v=h_psKQV-pU4}
}

@misc{GowersYouTubeInverseTheoremsAdditiveCombinatorics,
  author       = {Gowers, Timothy},
  title        = {Inverse theorems in additive combinatorics},
  howpublished = {YouTube video (KSoM Channel)},
  note         = {Video ID: Z3-v0apV\_t0. YouTube lists it as "Streamed 5 years ago" (accessed 2026-03-04).},
  url          = {https://www.youtube.com/watch?v=Z3-v0apV_t0}
}

@article{FuGuthMaldague2023,
  author  = {Fu, Yuqiu and Guth, Larry and Maldague, Dominique},
  title   = {Sharp superlevel set estimates for small cap decouplings of the parabola},
  journal = {Revista Matem{\'a}tica Iberoamericana},
  volume  = {39},
  number  = {3},
  pages   = {975--1004},
  year    = {2023},
  doi     = {10.4171/RMI/1393}
}

@article{Guth2016,
  author  = {Guth, Larry},
  title   = {A restriction estimate using polynomial partitioning},
  journal = {Journal of the American Mathematical Society},
  volume  = {29},
  number  = {2},
  pages   = {371--413},
  year    = {2016},
  doi     = {10.1090/jams/827}
}

@article{CarberyIliopoulouWang2024,
  author  = {Carbery, Anthony and Iliopoulou, Marina and Wang, Hong},
  title   = {Some sharp inequalities of {Mizohata--Takeuchi}-type},
  journal = {Revista Matemática Iberoamericana},
  volume  = {40},
  number  = {4},
  pages   = {1387--1418},
  year    = {2024},
  doi     = {10.4171/RMI/1463},
  eprint  = {2302.11877},
  archivePrefix = {arXiv}
}

@article{GuthHickmanIliopoulou2019,
    author  = {Guth, Larry and Hickman, Jonathan and Iliopoulou, Marina},
    title   = {Sharp estimates for oscillatory integral operators via polynomial partitioning},
    journal = {Acta Mathematica},
    volume  = {223},
    number  = {2},
    pages   = {251--376},
    year    = {2019},
    doi     = {10.4310/ACTA.2019.v223.n2.a2},
    eprint  = {1710.10349},
    archivePrefix = {arXiv},
    primaryClass  = {math.CA}
}

@misc{vanHandel2016,
    author       = {van Handel, Ramon},
    title        = {Probability in High Dimension},
    year         = {2016},
    howpublished = {Lecture notes for APC 550, Princeton University},
    url          = {https://web.math.princeton.edu/~rvan/APC550.pdf}
}

@misc{Guth2017Decoupling,
    author       = {Larry Guth},
    title        = {Topics in Analysis: Decoupling},
    year         = {2017},
    note         = {Lecture notes for MIT 18.118, Fall 2017},
    howpublished = {Massachusetts Institute of Technology},
    url          = {https://math.mit.edu/~lguth/Math118.html}
}

@book{Wolff2003,
    author    = {Thomas H. Wolff},
    title     = {Lectures on Harmonic Analysis},
    editor    = {Izabella {\L}aba and Carol Shubin},
    series    = {University Lecture Series},
    volume    = {29},
    publisher = {American Mathematical Society},
    address   = {Providence, RI},
    year      = {2003},
    isbn      = {978-0-8218-3449-7},
    doi       = {10.1090/ulect/029}
}
\end{document}